# Construction and properties of fuzzy-valued fractal interpolation function by using iterated function system

CholHui Yun[a*]   Hyang Choe[a]   MiGyong Ri[a]

[a]*Faculty of Mathematics, **Kim Il Sung** University, Pyongyang, Democratic People's Republic of Korea*

*ch.yun@ryongnamsan.edu.kp

In nature, there are many phenomena with both irregularity and uncertainty. Therefore, a fuzzy-valued fractal interpolation is more useful for modeling them than fuzzy interpolation or fractal interpolation. We construct fractal interpolation functions that interpolate a given data set of fuzzy numbers using an iterated function system and study the Hölder continuity of the constructed function. Firstly, we construct an iterated function system (IFS) using the given data set and construct a fuzzy-valued fractal interpolation function whose graph is the attractor of the constructed IFS. Next, we find the relationship between the fuzzy-valued fractal interpolation function and the fractal interpolation function of the level set of the data set. Finally, we prove the Hölder continuity of the fuzzy-valued fractal interpolation function.

*Keyword*s: Fractal interpolation function (FIF); fuzzy interpolation; iterated function system (IFS); fuzzy number; fuzzy fractal interpolation; Hölder continuity.

## 1. Introduction

Fractal interpolation function (FIF) is an interpolation function whose graph is a fractal set. The FIF is very useful to model objects with highly complicated structures. In 1986, Barnsley[1] introduced the FIF on the basis of an iterated function system (IFS).

The IFS is a family of transformations defined on a complete metric space X. We use these transformations of IFS to define the Hutchinson-Barnsley (HB) operator on the space of compact subsets of X and the Read-Bajraktarevic (RB) operator on a space of interpolation functions. The fixed point of HB operator is called an attractor of IFS. The attractor has the self-similarity. Furthermore, the fixed point of RB operator is just the FIF defined by Barnsley and the graph of the FIF is the attractor of IFS. This means that we can use the FIFs to model phenomena with self-similarity better than the classical interpolation functions such as polynomials and spline. Therefore, in many papers, constructions and properties of different kinds of FIFs such as recurrent FIF, hidden variable FIF and hidden variable recurrent FIF have been studied and applied to many fields such as computer graphics, signal process, image process and approximation theory, etc.[2-5]

Since the data given from experiments, observations and measurements etc. are perturbed by various factors, they have the uncertainty. It is important to precisely model the data with uncertainty in data processing and analysis.

To represent the degree of uncertainty in a mathematical way, in 1965, Zadeh[6] introduced the concept of fuzzy set. Furthermore, he proposed how to interpolate the data set of fuzzy numbers (Lowen (1990)[7]). To solve the problem, Lowen (1990)[7] constructed a fuzzy Lagrange interpolation polynomial which is continuous with respect to the Hausdorff metric between endographs of fuzzy sets. Kaleva (1994)[8] studied methods of calculating fuzzy Lagrange interpolation polynomial and fuzzy spline. Goghary (2005)[9] introduced a numerical method for calculating a fuzzy Hermite polynomial. Abbasbandy et al. (2008)[10] studied interpolation of the fuzzy data by a new family of spline functions.

Since the data obtained from highly irregular phenomena are perturbed by many factors, they have not only fractal properties but also uncertainty. Therefore, the interpolation which can represent both fractal characteristics and uncertainty is more appropriate for modeling them than fuzzy interpolation.[7-9]

A lot of researchers have studied the problem of applying the concept of fuzzification introduced by Zadeh[6] to fractals.[3, 11-17] Xiao et al [16, 17] proposed a method for modeling uncertain data using FIFs. They convert points of a given data into fuzzy numbers and then generated a family of FIFs using sets of start and end points of level sets of the fuzzy numbers. However, the family is not fuzzy-valued function with self-similarity.

In order to better describe objects with both self-similarity and uncertainty, we are going to investigate how to construct fuzzy-valued interpolation functions using IFS and their properties.

The remainder of this paper is organized as follows: in Section 2, we describe some preliminaries necessary to understand the paper. Section 3 describes the construction of fuzzy-valued fractal interpolation functions using iterated function systems; Section 4 studies the relationship between fuzzy-valued fractal interpolation functions and fractal interpolation functions according to level sets; Section 5 proves the Hölder continuity of fuzzy-valued fractal interpolation functions.

## 2. Preliminaries

Fuzzy number is a fuzzy set that is normal, upper semi-continuous, convex and compactly supported. Let us denote the space of fuzzy numbers over $\mathbf{R}$ by $\mathbf{R}_F$.

If $u \in \mathbf{R}_F$, then for any $\lambda \in [0,1]$, the $\lambda$-level set, $[u]^\lambda$ is a closed interval, which is denoted by $[u^-(\lambda), u^+(\lambda)]$. We define $d_\infty : \mathbf{R}_F \times \mathbf{R}_F \to \mathbf{R}$ as follows:

$$d_\infty(u,v) = \sup_{\lambda \in [0,1]} \max\{|u^-(\lambda) - v^-(\lambda)|, |u^+(\lambda) - v^+(\lambda)|\}, \quad u, v \in \mathbf{R}_F.$$

Then $d_\infty$ is a metric on $\mathbf{R}_F$ and $(\mathbf{R}_F, d_\infty)$ is a complete metric space.[18] $d_\infty$ is called a supreme metric on $\mathbf{R}_F$.

For any $u, v \in \mathbf{R}_F$ and $\lambda \in \mathbf{R}$, $u \oplus v, \lambda u, u \vee_g v$ are defined as follows by using their level sets:

$$[u \oplus v]^\lambda := [u]^\lambda + [v]^\lambda,$$
$$[\lambda u]^\lambda := \lambda [u]^\lambda,$$
$$[u \vee_g v]^\lambda := [\inf_{\beta \geq \lambda} \min\{u^-(\beta) - v^-(\beta), u^+(\beta) - v^+(\beta)\}, \sup_{\beta \geq \lambda} \max\{u^-(\beta) - v^-(\beta), u^+(\beta) - v^+(\beta)\}]$$

Let us denote the set of all fuzzy-valued continuous functions on a compact subset $\mathbf{K}$ of $\mathbf{R}$ by $C(\mathbf{K}, \mathbf{R}_F)$. We define $D : C(\mathbf{K}, \mathbf{R}_F) \times C(\mathbf{K}, \mathbf{R}_F) \to \mathbf{R}$ as follows:

$$D(f, g) = \sup_{x \in \mathbf{K}} d_\infty(f(x), g(x)), \quad f, g \in C(\mathbf{K}, \mathbf{R}_F),$$

where $d_\infty$ is the supreme metric on $\mathbf{R}_F$, i.e.

$$d_\infty(f(x), g(x)) = \sup_{\lambda \in [0,1]} \max\{|[f(x)]^+(\lambda) - [g(x)]^+(\lambda)|, |[f(x)]^-(\lambda) - [g(x)]^-(\lambda)|\}.$$

**Theorem 1.**[18]. $(C(\mathbf{K}, \mathbf{R}_F), D)$ *is a complete metric space.*

For $u \in C(\mathbf{K}, \mathbf{R}_F)$ and any $\lambda \in [0,1]$, a set
$$\{(x, (u(x))^-(\lambda)) \in \mathbf{R}^2 \mid x \in \mathbf{K}\} \cup \{(x, (u(x))^+(\lambda)) \in \mathbf{R}^2 \mid x \in \mathbf{K}\}$$
is called a $\lambda$-level set of $u$.

**Definition 1.**[5] Let $f : [a, b] \to \mathbf{R}$ be a function. If there are positive number $L$ and $0 \leq \lambda \leq 1$ such that
$$\forall x, y \in [a, b], |f(x) - f(y)| \leq L|x - y|^\lambda,$$
then $f$ is called a Hölder continuous function with exponent $\lambda$. $L$ is called a Hölder coefficient of $f$ and denoted by $H_f$. $\lambda$ is called a Hölder exponent of $f$ and denoted by $\lambda_f$. If $\lambda_f = 1$, then $f$ is called a Lipschtz function and $H_f$ is denoted by $L_f$. Especially, if $L_f \in [0,1)$, then $f$ is called a contraction and $L_f$ is denoted by $c_f$.

## 3. Construction of fuzzy-valued fractal interpolation function.

In this section, we consider a method of constructing fuzzy-valued interpolation function using IFS.

A data set $P$ of fuzzy numbers is given as follows:
$$P = \{(x_i, u_i) \in \mathbf{R} \times \mathbf{R_F} \mid i = 0, 1, \cdots, n\}, (x_0 < x_1 < \cdots x_n).$$
Let us denote $I = [x_0, x_n]$, $I_i = [x_{i-1}, x_i]$ for $i = 1, \cdots, n$ and $l_i : I \to I_i$, $i = 1, \cdots, n$ be contractive homeomorphisms such that
$$l_i(x_0) = x_{i-1}, l_i(x_n) = x_i. \quad (1)$$
We define mappings $F_i : I \times \mathbf{R_F} \to \mathbf{R_F}$, $i = 1, \cdots, n$ as follows:
$$(x, u) \in I \times \mathbf{R_F}, \quad F_i(x, u) = s_i u \oplus q_i(l_i(x)), \quad (2)$$
where $s_i$ is a constant such that $0 \le s_i < 1$ and $q_i : I_i \to \mathbf{R_F}$ is a Lipschitz mapping satisfying the following condition:
$$F_i(x_0, u_0) = u_{i-1}, F_i(x_n, u_n) = u_i. \quad (3)$$
The following is an example of $q_i(x)$.
$$q_i(x) = b_i(x) \vee_g s_i g_i(l_i^{-1}(x)),$$
$$g_i(x) = \frac{x - x_0}{x_n - x_0} u_n \oplus \frac{x - x_n}{x_0 - x_n} u_0, \quad b_i(x) = \frac{x - x_{i-1}}{x_i - x_{i-1}} u_i \oplus \frac{x - x_i}{x_{i-1} - x_i} u_{i-1}.$$
In this case, $F_i(x, u) = s_i (u \vee_g g_i(x)) \oplus b_i(l_i(x))$.

We define a metric $d_{\max} : (I \times \mathbf{R_F}) \times (I \times \mathbf{R_F}) \to \mathbf{R}$ as follows:
$$d_{\max}((x, u), (y, v)) = \max\{|x - y|, d_\infty(u, v)\}, \quad (x, u), (y, v) \in I \times \mathbf{R_F}.$$
Then $(I \times \mathbf{R_F}, d_{\max})$ is a complete metric space since it is a product of complete metric spaces.

We define transformations $w_i : I \times \mathbf{R_F} \to I_i \times \mathbf{R_F}$, $i = 1, \cdots, n$ by
$$w_i(x, u) = (l_i(x), F_i(x, u)), \quad (x, u) \in I \times \mathbf{R_F}.$$

**Theorem 2.** *There exists a metric $d_\theta$ equivalent to $d_{\max}$ such that for $i = 1, \cdots, n$, $w_i$ are contractive with respect to $d_\theta$.*

**Proof.** We define a metric $d_\theta : (I \times \mathbf{R_F}) \times (I \times \mathbf{R_F}) \to \mathbf{R}$ by
$$d_\theta((x, u), (x', u')) = |x - x'| + \theta d_\infty(u, u'),$$
where $\theta$ is a real number satisfying $0 < \theta < \max_{i=1,\cdots,n}\{(1 - c_{l_i})/(L_{q_i} \cdot c_{l_i})\}$. Then the metric $d_\theta$ is equivalent to $d_{\max}$.

In fact, in the case of $0 < \theta < 1$, we have
$$\theta d_{\max}((x_1, u_1), (x_2, u_2)) = \theta \max\{|x_1 - x_2|, d_\infty(u_1, u_2)\} \le \theta |x_1 - x_2| + \theta d_\infty(u_1, u_2)$$
$$\le d_\theta((x_1, u_1), (x_2, u_2))$$
$$= |x_1 - x_2| + \theta d_\infty(u_1, u_2) \le |x_1 - x_2| + d_\infty(u_1, u_2)$$
$$\le 2 \max\{|x_1 - x_2|, d_\infty(u_1, u_2)\} = 2 d_{\max}((x_1, u_1), (x_2, u_2))$$

and in the case of $\theta \ge 1$, we get
$$d_{\max}((x_1, u_1), (x_2, u_2)) = \max\{|x_1 - x_2|, d_\infty(u_1, u_2)\}$$
$$\le d_\theta((x_1, u_1), (x_2, u_2))$$
$$= |x_1 - x_2| + \theta d_\infty(u_1, u_2) \le 2 \max\{\theta |x_1 - x_2|, \theta d_\infty(u_1, u_2)\}$$
$$= 2\theta d_{\max}((x_1, u_1), (x_2, u_2)).$$

For $(x, u), (x', u') \in I \times \mathbf{R_F}$, we obtain
$$d_\theta(w_i(x, u), w_i(x', u')) = d_\theta((l_i(x), F_i(x, u)), (l_i(x'), F_i(x', u')))$$
$$= |l_i(x) - l_i(x')| + \theta d_\infty(F_i(x, u), F_i(x', u'))$$
$$\le c_{l_i} |x - x'| + \theta d_\infty(s_i u \oplus q_i(l_i(x)), s_i u' \oplus q_i(l_i(x')))$$
$$\le c_{l_i} |x - x'| + \theta(d_\infty(s_i u, s_i u') + d_\infty(q_i(l_i(x)), q_i(l_i(x'))))$$

$$\leq c_{l_i} |x-x'| + \theta |s_i| d_\infty(u,u') + \theta L_{q_i} c_{l_i} |x-x'|$$
$$= (c_{l_i} + \theta L_{q_i} c_{l_i}) |x-x'| + \theta \cdot s_i \cdot d_\infty(u,u')$$
$$\leq c_{w_i} (|x-x'| + \theta d_\infty(u,u'))$$
$$= c_{w_i} d_\theta((x,u),(x',u')),$$

where $c_{w_i} = \max\{c_{l_i} + \theta L_{q_i} c_{l_i}, s_i\} < 1$. Therefore, for $i = 1,\cdots,n$, $w_i$ are contractive. □

Since $(I \times \mathbf{R_F}, d_{\max})$ is a complete metric space, this theorem shows that $(I \times \mathbf{R_F}, d_\theta)$ is also complete.

Hence, by Theorem 2, $\{I \times \mathbf{R_F}; w_i = (l_i, F_i), i = 1,\cdots,n\}$ is a hyperbolic iterated function system (IFS). Now we denote the family of all compact subsets of $I \times \mathbf{R_F}$ by $H(I \times \mathbf{R_F})$, where the distance is the Hausdorff distance. The Hutchinson-Barnsley operator $W : H(I \times \mathbf{R_F}) \to H(I \times \mathbf{R_F})$ is defined as follows:

$$W(A) = \bigcup_{i=1}^{n} w_i(A), \quad A \in H(I \times \mathbf{R_F}).$$

Since for $i = 1,\cdots,n$, $w_i$ are contractive, $W$ is also contractive. Therefore, $W$ has a unique fixed point $B \in H(I \times \mathbf{R_F})$, which is the attractor of the IFS.

**Theorem 3.** *The attractor $B$ of the IFS $\{I \times \mathbf{R_F}; w_i = (l_i, F_i), i = 1,\cdots,n\}$ is the graph of a fuzzy-valued continuous function.*

**Proof.** Let us define a space of fuzzy-valued continuous functions on $I$ by
$$C^*(I, \mathbf{R_F}) = \{\varphi : I \to \mathbf{R_F} \mid \varphi \text{ is contiuous and } \varphi(x_i) = u_i, i = 0, 1, \cdots, n\}.$$

Then $C^*(I, \mathbf{R_F}) \subset C(I, \mathbf{R_F})$ and by Theorem 1, $(C^*, D)$ is a complete metric space. In fact, if $\{f_n\}$ is a Cauchy sequence in $(C^*, D)$, then because of $\{f_n\} \subset C^* \subset C$ and Theorem 1, there exists a $f \in C(I, \mathbf{R_F})$ such that $f_n \to f$ as $n \to \infty$. Since $\{f_n\} \subset C^*$, we have $f_n(x_i) = u_i$, $i = 0, 1, \cdots, n$ and $f(x_i) = \lim_{n \to \infty} f_n(x_i) = \lim_{n \to \infty} u_i = u_i$, that is, $f \in C^*(I, \mathbf{R_F})$. Therefore, $\{f_n\}$ is a convergent sequence in $(C^*, D)$.

Now, for $\varphi \in C^*(I, \mathbf{R_F})$, we define a function $T(\varphi) : I \to \mathbf{R_F}$ by
$$(T(\varphi))(x) = F_i(l_i^{-1}(x), \varphi(l_i^{-1}(x))), \quad x \in I_i.$$

Then $T(\varphi) \in C^*(I, \mathbf{R_F})$. In fact, from the definition of $F_i$, $T(\varphi)$ is continuous on $I$ and we have
$$(T(\varphi))(x_i) = (T(\varphi))(l_i(x_n)) = F_i(x_n, \varphi(x_n)) = F_i(x_n, u_n) = u_i.$$

Therefore, RB operator $T : C^*(I, \mathbf{R_F}) \to C^*(I, \mathbf{R_F})$ is well defined. Moreover, $T$ is contractive. In fact, for any $\varphi_1, \varphi_2 \in C^*(I, \mathbf{R_F})$, we have
$$D(T(\varphi_1), T(\varphi_2)) = \sup_{x \in I} d_\infty((T(\varphi_1))(x), (T(\varphi_2))(x))$$

and for each $x \in I_i$, we get
$$d_\infty((T\varphi_1)(x), (T\varphi_2)(x)) = d_\infty(F_i(l_i^{-1}(x), \varphi_1(l_i^{-1}(x))), F_i(l_i^{-1}(x), \varphi_2(l_i^{-1}(x))))$$
$$= d_\infty(s_i \varphi_1(l_i^{-1}(x)) \oplus q_i(x), s_i \varphi_2(l_i^{-1}(x)) \oplus q_i(x))$$
$$= d_\infty(s_i \varphi_1(l_i^{-1}(x)), s_i \varphi_2(l_i^{-1}(x)))$$
$$= s_i d_\infty(\varphi_1(l_i^{-1}(x)), \varphi_2(l_i^{-1}(x)))$$
$$\leq s D(\varphi_1, \varphi_2),$$

where $s = \max_{i=1,\cdots,n}\{s_i\}(<1)$. Hence, we obtain

$$D(T(\varphi_1), T(\varphi_2)) = \sup_{x\in I} d_\infty((T(\varphi_1))(x), (T(\varphi_2))(x)) \le sD(\varphi_1, \varphi_2),$$

that is, $T$ is a contractive operator. Then by Banach fixed point theorem, there exists a unique fixed point $f(\in C^*(I, \mathbf{R_F}))$ of $T$, i.e., we have

$$f(x) = Tf(x) = F_i(l_i^{-1}(x), f(l_i^{-1}(x))) = s_i f(l_i^{-1}(x)) \oplus q_i(x), \ x \in I_i. \quad (4)$$

Then $Gr\ f = B$, where $Gr\ f = \{(x, f(x)) \in I \times \mathbf{R_F} : x \in I\}$. In fact,

$$\begin{aligned}
Gr\ f\big|_{I_i} &= \{(x, f(x)) \in I \times \mathbf{R_F} : x \in I_i\} \\
&= \{(x, F_i(l_i^{-1}(x), f(l_i^{-1}(x)))) \in I \times \mathbf{R_F} : x \in I_i\} \\
&= \{(l_i(x), F_i(x, f(x))) \in I \times \mathbf{R_F} : x \in I\} \\
&= \{w_i(x, f(x)) : x \in I\} \\
&= w_i(Gr\ f).
\end{aligned}$$

Therefore, $Gr\ f = \bigcup_{i=1}^{n} w_i(Gr\ f)$, that is, $Gr\ f$ is the attractor of the IFS $\{I \times \mathbf{R_F}; w_i = (l_i, F_i), i = 1,\cdots,n\}$. Since the attractor of the IFS is unique, $Gr\ f = B$. □

**Example 1**. Let a data set $P$ of fuzzy numbers be given as follows:

$$P = \{(x_i, u_i) \in \mathbf{R} \times \mathbf{R_F} \mid i = 0, 1, \cdots, 4\} \left(x_i = \frac{i}{4}\right),$$

$$u_0(y) = \begin{cases} y, & 0 \le y \le 1 \\ 2-y, & 1 < y \le 2 \\ 0, & \text{otherwise} \end{cases}, \quad u_1(y) = \begin{cases} 1 - 0.25(y-3)^3, & 1 \le y \le 5 \\ 0, & \text{otherwise} \end{cases},$$

$$u_2(y) = \begin{cases} e^{-(y-3)^2}, & 0 \le y \le 6 \\ 0, & \text{otherwise} \end{cases}, \quad u_3(y) = \begin{cases} 2y-4, & 2 \le y \le 2.5 \\ 6-2y, & 2.5 < y \le 3 \\ 0, & \text{otherwise} \end{cases},$$

$$u_4(y) = \begin{cases} (y-1)^2/9, & 1 \le y \le 4 \\ (y-7)^2/9, & 4 < y \le 7 \\ 0, & \text{otherwise} \end{cases}.$$

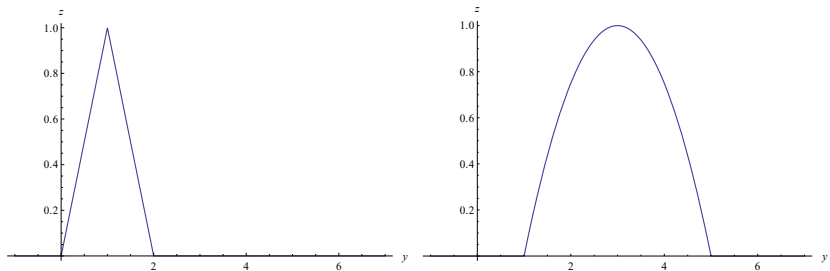

1)                        2)

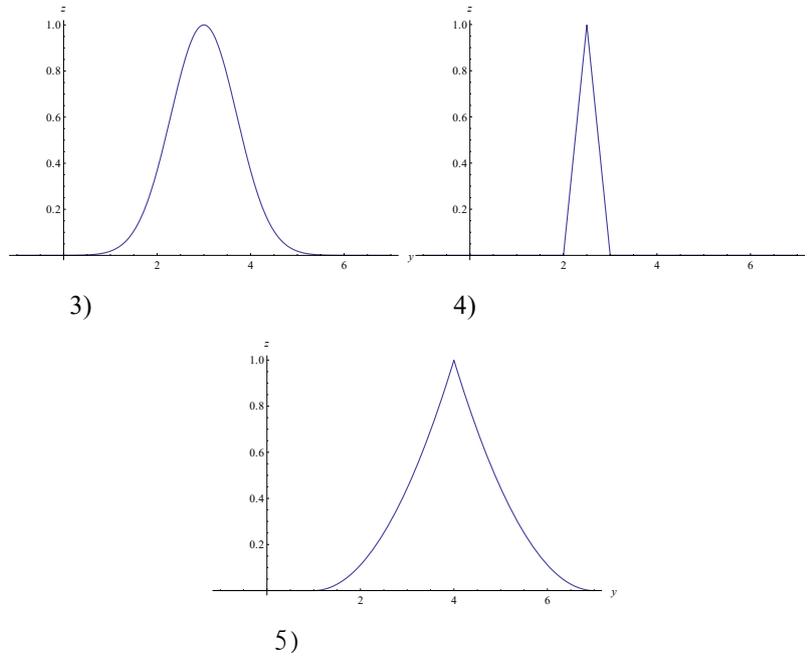

Fig. 1. 1), 2), 3), 4), 5) are the graphs of $u_0(y)$, $u_1(y)$, $u_2(y)$, $u_3(y)$, $u_4(y)$.

The vertical scaling factors $s_i$ are as follows:

$$s_1 = 0.3, \ s_2 = 0.7, \ s_3 = 0.4, \ s_4 = 0.8$$

Moreover, Lipschitz mappings $q_i : I_i \to \mathbf{R_F}$, $i = 1, 2, 3, 4$ satisfying (3) are given as follows:

$$q_i(x) = \left( \frac{x - x_{i-1}}{x_i - x_{i-1}} u_i \oplus \frac{x - x_i}{x_{i-1} - x_i} u_{i-1} \right) \mathbf{v}_g \left[ s_i \left( \frac{l_i^{-1}(x) - x_0}{x_n - x_0} u_n \oplus \frac{l_i^{-1}(x) - x_n}{x_0 - x_n} u_0 \right) \right].$$

Figure 2 shows the graph of the fuzzy-valued fractal interpolation function $f$.

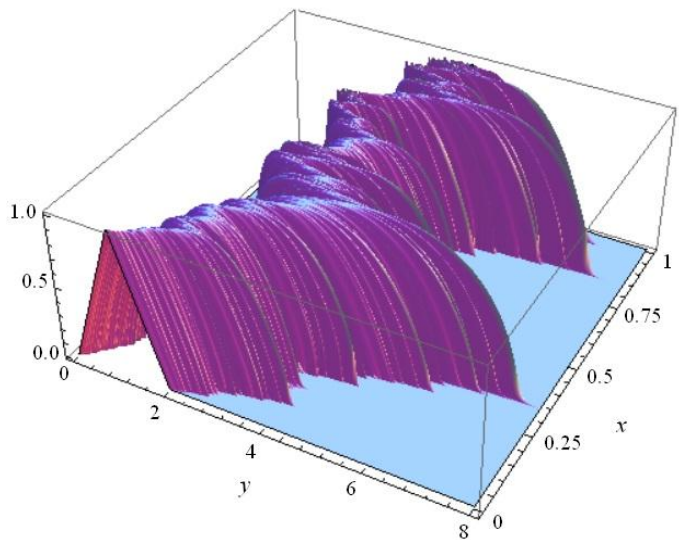

Fig. 2. The graph of the fuzzy-valued fractal interpolation function $f$

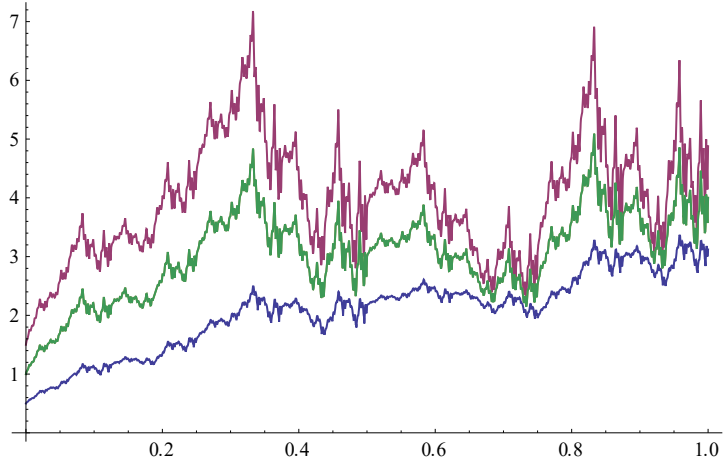

Fig. 3. The graphs of 1-level set and 0.5-level set

(Blue- $(f(x))^-(0.5)$, Purple- $(f(x))^+(0.5)$, Green- $(f(x))^+(1) = (f(x))^-(1)$)

Figure 3 shows the graphs of 1-level set and 0.5-level set of $f$.

## 4. Relationship between fuzzy-valued fractal interpolation function and fractal interpolation function interpolating the level sets.

In this section, we consider the relationship between the level sets of fuzzy-valued fractal interpolation function $f$ and the fractal interpolation functions interpolating the start and end points of the level sets of the data set $P$ of fuzzy numbers. For all $\lambda \in [0, 1]$, let us denote as follows:

$$P_\lambda^- = \{(x_i, u_i^-(\lambda)) \in \mathbf{R}^2 \mid i = 0, 1, \cdots, n\}, \quad P_\lambda^+ = \{(x_i, u_i^+(\lambda)) \in \mathbf{R}^2 \mid i = 0, 1, \cdots, n\}.$$

We define the mappings $(F_i)_\lambda^- : I \times \mathbf{R} \to \mathbf{R}$ and $(F_i)_\lambda^+ : I \times \mathbf{R} \to \mathbf{R}$, $i = 1, \cdots, n$ as follows:

$$(F_i)_\lambda^-(x, y) = s_i y + (q_i(l_i(x)))^-(\lambda), \quad (F_i)_\lambda^+(x, y) = s_i y + (q_i(l_i(x)))^+(\lambda),$$

where $s_i$, $l_i$ and $q_i$ are given by (1) and (2).

From (2), we have

$$\begin{aligned}[][F_i(x_0, u_0)]^\lambda &= s_i[u_0]^\lambda + [q_i(l_i(x_0))]^\lambda \\
&= s_i[u_0^-(\lambda), u_0^+(\lambda)] + [(q_i(l_i(x)))^-(\lambda), (q_i(l_i(x)))^+(\lambda)] \\
&= [s_i u_0^-(\lambda) + (q_i(l_i(x)))^-(\lambda), s_i u_0^+(\lambda) + (q_i(l_i(x)))^+(\lambda)] \\
&= [(F_i)_\lambda^-(x_0, u_0^-(\lambda)), (F_i)_\lambda^+(x_0, u_0^+(\lambda))]\end{aligned}$$

and from (3), we get

$$[F_i(x_0, u_0)]^\lambda = [u_{i-1}]^\lambda = [u_{i-1}^-(\lambda), u_{i-1}^+(\lambda)].$$

Hence, we obtain

$$(F_i)_\lambda^-(x_0, u_0^-(\lambda)) = u_{i-1}^-(\lambda), \quad (F_i)_\lambda^+(x_0, u_0^+(\lambda)) = u_{i-1}^+(\lambda).$$

Similarly, we get

$$(F_i)_\lambda^-(x_n, u_n^-(\lambda)) = u_i^-(\lambda), \quad (F_i)_\lambda^+(x_n, u_n^+(\lambda)) = u_i^+(\lambda).$$

Then $\{\mathbf{R}^2; (l_i, (F_i)_\lambda^-), i = 1, \cdots, n\}$ and $\{\mathbf{R}^2; (l_i, (F_i)_\lambda^+), i = 1, \cdots, n\}$ are hyperbolic iterated function systems and the attractors of these IFSs are the graphs of continuous functions $g_\lambda^-$ and $g_\lambda^+$ interpolating the data sets $P_\lambda^-$ and $P_\lambda^+$, respectively:

$$g_\lambda^-(x) = s_i g_\lambda^-(l_i^{-1}(x)) + (q_i(x))^-(\lambda), \quad g_\lambda^+(x) = s_i g_\lambda^+(l_i^{-1}(x)) + (q_i(x))^+(\lambda).$$

These functions are the fixed points of the RB operators $T_\lambda^-$, $T_\lambda^+$ defined on complete metric spaces $\{\varphi \in C(I, \mathbf{R}) : \varphi(x_i) = u_i^-(\lambda), i = 0, 1, \cdots, n\}$, $\{\psi \in C(I, \mathbf{R}) : \psi(x_i) = u_i^+(\lambda), i = 0, 1, \cdots, n\}$, respectively.

$$T_\lambda^-(\varphi) = s_i \varphi(l_i^{-1}(x)) + (q_i(x))^-(\lambda), \tag{5}$$

$$T_\lambda^+(\psi) = s_i \psi(l_i^{-1}(x)) + (q_i(x))^+(\lambda). \tag{6}$$

On the other hand, by (4), the fuzzy-valued fractal interpolation function $f$ interpolating $P$ satisfies the following:

$$[f(x)]^\lambda = [(f(x))^-(\lambda), (f(x))^+(\lambda)]$$
$$= s_i [f(l_i^{-1}(x))]^\lambda + [q_i(x)]^\lambda$$
$$= [s_i (f(l_i^{-1}(x)))^-(\lambda) + (q_i(x))^-(\lambda), \; s_i (f(l_i^{-1}(x)))^+(\lambda) + (q_i(x))^+(\lambda)],$$

i.e.

$$(f(x))^-(\lambda) = s_i (f(l_i^{-1}(x)))^-(\lambda) + (q_i(x))^-(\lambda) \tag{7}$$

$$(f(x))^+(\lambda) = s_i (f(l_i^{-1}(x)))^+(\lambda) + (q_i(x))^+(\lambda). \tag{8}$$

Because of (5)-(8) and uniqueness of fixed point of a contractive mapping on a complete metric space, the following holds:

$$(f(x))^-(\lambda) = g_\lambda^-(x), \; (f(x))^+(\lambda) = g_\lambda^+(x).$$

Therefore, we obtain the following theorem.

**Theorem 4.** *Let $f$ be a fuzzy-valued fractal interpolation function interpolating the fuzzy-valued data set $P$. Then for each $\lambda \in [0, 1]$, $(f(\cdot))^-(\lambda)$ and $(f(\cdot))^+(\lambda)$ are the fractal interpolation functions interpolating $P_\lambda^-$ and $P_\lambda^+$, respectively.*

**Example 2.** Figure 4 shows the graphs of fractal interpolation functions $g_{0.5}^-$, $g_{0.5}^+$, $g_1^- = g_1^+$ that interpolate 0.5-level set and 1-level set

$$P_{0.5}^- = \{(x_i, u_i^-(0.5)) \in \mathbf{R}^2 \mid i = 0, 1, \cdots, 4\}, \; P_{0.5}^+ = \{(x_i, u_i^+(0.5)) \in \mathbf{R}^2 \mid i = 0, 1, \cdots, 4\},$$
$$P_1^- = P_1^+ = \{(x_i, u_i^+(1)) \in \mathbf{R}^2 \mid i = 0, 1, \cdots, 4\}$$

of the data set of fuzzy numbers given in the Example 1.

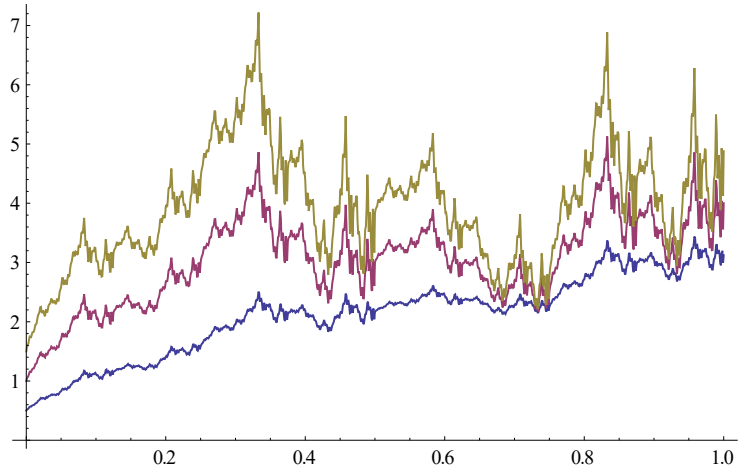

Fig. 4. The graphs of $g_{0.5}^-$, $g_{0.5}^+$, $g_1^- = g_1^+$

(Blue- $g_{0.5}^-$, Green- $g_{0.5}^+$, Purple- $g_1^- = g_1^+$ )

Figure 3 and figure 4 show that the graphs of $(f(\cdot))^-(0.5)$, $(f(\cdot))^+(0.5)$, $(f(\cdot))^-(1) = (f(\cdot))^+(1)$ are same as the graphs of $g_{0.5}^-$, $g_{0.5}^+$, $g_1^- = g_1^+$, respectively.

## 5. Hölder continuity of the fuzzy-valued fractal interpolation function

In this section, we consider a Hölder continuity of the fuzzy-valued fractal interpolation function.

First of all, we prove the boundedness and Hölder continuity of the fractal interpolation function. Let us denote the fractal interpolation function interpolating the data set $\{(x_i, y_i): i = 0, 1, \cdots, n\}$ by $f_\mathbf{y}$:

$$f_\mathbf{y}(x) = s_i f_\mathbf{y}(l_i^{-1}(x)) + q_i(x), \quad x \in I_i,$$

where $\mathbf{y} = (y_0, y_1, \cdots, y_n)$, $|s_i| < 1$ and $q_i(x): I_i \to \mathbf{R}$ are Lipschitz functions satisfying the following:

1) There is a constant $\rho > 0$ such that $L_{q_i} < \rho, i = 1, \cdots, n$

2) $q_i(x_{i-1}) = y_{i-1} - s_i y_0$, $q_i(x_i) = y_i - s_i y_N$ \hfill (9)

and $l_i$ are similitude, i.e., there exists $c_{l_i} \in [0, 1)$ such that

$$\forall x, y \in I, \ |l_i(x) - l_i(y)| = c_{l_i} |x - y|.$$

Let us denote as follows:

$$c_{\min} = \min_{i=1,\cdots,n} \{c_{l_i}\}, \ c_{\max} = \max_{i=1,\cdots,n} \{c_{l_i}\}, \ L_q = \max_{i=1,\cdots,n} \{L_{q_i}\}.$$

**Theorem 5.** *If there exists $A \in \mathbf{R}$ such that $|y_i| \leq A$ for $i = 1, \cdots, n$, then the fractal interpolation function $f_\mathbf{y}$ has the following properties:*

1) *There exists a constant $\alpha$ independent of $\mathbf{y}$ such that $\|f_\mathbf{y}\|_\infty \leq \alpha$.*

2) *$f_\mathbf{y}$ is a Hölder continuous function, i.e.*

$$\exists K > 0, \ \tau(0 < \tau \leq 1): \ \forall x, x' \in I, \ |f_\mathbf{y}(x) - f_\mathbf{y}(x')| \leq K |x - x'|^\tau,$$

*where $K$ and $\tau$ are the constants independent of $\mathbf{y}$.*

**Proof.** 1) Since for any $x \in I_i$, there exists $x' \in I$ such that

$$f_\mathbf{y}(x) = s_i f_\mathbf{y}(x') + q_i(x),$$

we get

$$|f_\mathbf{y}(x)| \leq s_i |f_\mathbf{y}(x')| + |q_i(x)| \leq s \|f_\mathbf{y}\|_\infty + \|q_i\|_\infty.$$

Moreover, $q_i$ is a Lipschitz function satisfying (9). Therefore, for any $x \in I_i$, we have

$$|q_i(x) - q_i(x_{i-1})| \leq L_{q_i} |x - x_{i-1}| \leq \rho \cdot c_{l_i}(x_n - x_0) \leq \rho \cdot c_{\max}(x_n - x_0),$$
$$|q_i(x)| \leq \rho \cdot c_{\max}(x_n - x_0) + |q_i(x_{i-1})| \leq \rho \cdot c_{\max}(x_n - x_0) + (1+s)A,$$
$$\|q_i\|_\infty \leq \rho \cdot c_{\max}(x_n - x_0) + (1+s)A.$$

Hence, we obtain

$$|f_\mathbf{y}(x)| \leq s \|f_\mathbf{y}\|_\infty + \rho \cdot c_{\max}(x_n - x_0) + (1+s)A,$$
$$\|f_\mathbf{y}\|_\infty \leq s \|f_\mathbf{y}\|_\infty + \rho \cdot c_{\max}(x_n - x_0) + (1+s)A,$$
$$\|f_\mathbf{y}\|_\infty \leq \alpha,$$

where $\alpha = \dfrac{\rho \cdot c_{\max}(x_n - x_0) + (1+s)A}{1-s}$.

2) Let us denote $I_{r_1 r_2 \cdots r_m} = l_{r_m} \circ l_{r_{m-1}} \circ \cdots \circ l_{r_1}(I), \ r_i \in \{1, 2, \cdots, n\}$.

For any $x, x' \in I \ (x < x')$, firstly, we consider the case where $x, x' \in I_i$. In this case, there exist $n', n'' \ (n' \neq n'')$ such that $x \in I_{t_{n'} t_{n-1} \cdots t_1}, \ x' \in I_{t_{n''} t_{n-1} \cdots t_1}$.

Let us denote $I_{t_{n''} t_{n-1} \cdots t_1} =: [x'', x''']$. For an interval $[x, x'']$, there exists $m \ (n \leq m)$ such that $I_{r_1 r_2 \cdots r_m} \subset [x, x''] \subset I_{r_2 \cdots r_m}$, where $r_{m-n+2} = t_{n-1}, \cdots, r_m = t_1$.

We have

$$|f_\mathbf{y}(x) - f_\mathbf{y}(x'')| = |f_\mathbf{y}(l_{r_m}(l_{r_m}^{-1}(x))) - f_\mathbf{y}(l_{r_m}(l_{r_m}^{-1}(x'')))|$$

$$= |[s_{r_m} f_{\mathbf{y}}(l_{r_m}^{-1}(x)) + q_{r_m}(x)] - [s_{r_m} f_{\mathbf{y}}(l_{r_m}^{-1}(x'')) + q_{r_m}(x'')]|$$

$$\leq s \cdot |f_{\mathbf{y}}(l_{r_m}^{-1}(x)) - f_{\mathbf{y}}(l_{r_m}^{-1}(x''))| + L_q \cdot |x - x''|$$

$$\leq s \cdot [s \cdot |f_{\mathbf{y}}(l_{r_{m-1}r_m}^{-1}(x)) - f_{\mathbf{y}}(l_{r_{m-1}r_m}^{-1}(x''))| + L_q \cdot |l_{r_m}^{-1}(x) - l_{r_m}^{-1}(x'')|] + L_q \cdot |x - x''|$$

$$\leq s^2 \cdot |f_{\mathbf{y}}(l_{r_{m-1}r_m}^{-1}(x)) - f_{\mathbf{y}}(l_{r_{m-1}r_m}^{-1}(x''))| + L_q |x - x''| \left(1 + \frac{s}{c_{\min}}\right)$$

$$\leq \cdots \leq$$

$$\leq s^{m-1} \cdot |f_{\mathbf{y}}(l_{r_2 \cdots r_m}^{-1}(x)) - f_{\mathbf{y}}(l_{r_2 \cdots r_m}^{-1}(x''))| + L_q \cdot |x - x''| \left[1 + \frac{s}{c_{\min}} + \cdots + \left(\frac{s}{c_{\min}}\right)^{m-2}\right]$$

$$\leq 2 s^{m-1} \alpha + L_q \cdot |x - x''| \left[1 + \frac{s}{c_{\min}} + \cdots + \left(\frac{s}{c_{\min}}\right)^{m-2}\right]$$

Let us denote $\delta = s / c_{\min}$. Then we get

$$|f_{\mathbf{y}}(x) - f_{\mathbf{y}}(x'')| \leq 2\alpha \delta^{m-1} c_{\min}^{m-1} + L_q \cdot |x - x''| (1 + \delta + \cdots + \delta^{m-2}).$$

Since $I_{r_1 r_2 \cdots r_m} \subset [x, x''] \subset I_{r_2 \cdots r_m}$, we obtain

$$c_{\min}^m (x_n - x_0) \leq c_{l_1} \cdots c_{l_{r_m}} (x_n - x_0) \leq |x - x''| \leq c_{l_{r_2}} \cdots c_{l_{r_m}} (x_n - x_0) \leq c_{\max}^{m-1}(x_n - x_0)$$

and

$$\left| f_{\mathbf{y}}(x) - f_{\mathbf{y}}(x'') \right| \leq \frac{2\alpha \delta^{m-1}}{c_{\min}(x_n - x_0)} \cdot |x - x''| + \rho |x - x''| (1 + \delta + \cdots + \delta^{m-2})$$

$$\leq M |x - x''| (1 + \delta + \cdots + \delta^{m-1}),$$

where $M = \max\left\{\dfrac{2\alpha}{c_{\min}(x_n - x_0)}, \rho\right\}$.

Now we consider the following cases:

(1) In the case of $\delta < 1$, i.e. $s < c_{\min}$, we get

$$|f_{\mathbf{y}}(x) - f_{\mathbf{y}}(x'')| \leq Q |x - x''|^\tau,$$

where $Q = \dfrac{M}{1 - \delta}$ and $\tau = 1$.

(2) In the case of $\delta = 1$, i.e. $s = c_{\min}$, we obtain

$$|f_{\mathbf{y}}(x) - f_{\mathbf{y}}(x'')| \leq Mm |x - x''|.$$

Moreover, we have

$$|x - x''| \leq c_{\max}^{m-1}(x_n - x_0),$$

$$m - 1 \leq \frac{\ln[|x - x''| / (x_n - x_0)]}{\ln c_{\max}}$$

and since $0 < -x^\omega \ln x \leq \dfrac{1}{\omega \cdot e}$ for any $0 < x < 1$ and $\omega > 0$, for any $0 < \tau < 1$, we obtain

$$m |x - x''| \leq \left(1 + \frac{\ln[|x - x''| / (x_n - x_0)]}{\ln c_{\max}}\right) |x - x''|$$

$$\leq |x - x''| + \frac{-\left|\dfrac{x - x''}{x_n - x_0}\right|^{1-\tau} \ln \left|\dfrac{x - x''}{x_n - x_0}\right|}{|\ln c_{\max}|} |x - x''|^\tau (x_n - x_0)^{1-\tau}$$

$$\leq \left|\frac{x - x''}{x_n - x_0}\right|^{1-\tau} \cdot (x_n - x_0)^{1-\tau} |x - x''|^\tau + \frac{2^{1-\tau}(x_n - x_0)^{1-\tau}}{|\ln c_{\max}|(1-\tau)e} |x - x''|^\tau$$

$$\leq \left(1 + \frac{1}{|\ln c_{\max}|(1-\tau)e}\right)(x_n - x_0)^{1-\tau} |x - x''|^\tau$$

$$\leq \left(1 + \frac{1}{|\ln c_{\max}|(1-\tau)\mathrm{e}}\right)\max\{1,\ x_n - x_0\}|x - x''|^{\tau}$$

Therefore, we obtain
$$|f_{\mathbf{y}}(x) - f_{\mathbf{y}}(x'')| \leq Q|x - x''|^{\tau},$$
where $Q = M\left(1 - \dfrac{1}{(1-\tau)\mathrm{e}\ln c_{\max}}\right)\max\{1,\ x_n - x_0\}$.

(3) In the case of $\delta > 1$, i.e. $s > c_{\min}$, we get
$$|f_{\mathbf{y}}(x) - f_{\mathbf{y}}(x'')| \leq \frac{M\delta^m}{\delta - 1}|x - x''|.$$

Let us denote $\tau := \dfrac{\ln \delta}{\ln c_{\max}} + 1$. Then $\tau < 1$. Since $|x - x''| \leq c_{\max}^{m-1}(x_n - x_0)$, we have
$$\ln\left|\frac{x - x''}{x_n - x_0}\right| \leq (m-1)\ln c_{\max} = -\frac{(m-1)\ln \delta}{1-\tau},$$
$$(1-\tau)\ln\left|\frac{x - x''}{x_n - x_0}\right| \leq -(m-1)\ln \delta,$$
$$\ln\left|\frac{x - x''}{x_n - x_0}\right|^{1-\tau} \leq \ln\frac{1}{\delta^{m-1}},$$
$$\delta^{m-1}\left|\frac{x - x''}{x_n - x_0}\right|^{1-\tau} \leq 1,$$
$$\delta^{m-1}|x - x''| \leq (x_n - x_0)^{1-\tau}|x - x''|^{\tau} \leq \max\{1,\ x_n - x_0\}|x - x''|^{\tau}.$$

Hence, we get
$$|f_{\mathbf{y}}(x) - f_{\mathbf{y}}(x'')| \leq Q|x - x''|^{\tau},$$
where $Q = \dfrac{M\delta}{\delta - 1}\max\{1,\ x_n - x_0\}$.

Similarly, for $x''$ and $x'$, we have
$$|f_{\mathbf{y}}(x') - f_{\mathbf{y}}(x'')| \leq Q|x' - x''|^{\tau}.$$

Hence, we get
$$|f_{\mathbf{y}}(x) - f_{\mathbf{y}}(x')| \leq |f_{\mathbf{y}}(x) - f_{\mathbf{y}}(x'')| + |f_{\mathbf{y}}(x'') - f_{\mathbf{y}}(x')|$$
$$\leq Q|x - x''|^{\tau} + Q|x' - x''|^{\tau}$$
$$\leq 2Q|x - x'|^{\tau}$$

Next, in the case where $x \in I_i$, $x' \in I_j$ $(i < j)$, we obtain
$$|f_{\mathbf{y}}(x) - f_{\mathbf{y}}(x')| \leq |f_{\mathbf{y}}(x) - f_{\mathbf{y}}(x_i)| + |f_{\mathbf{y}}(x_i) - f_{\mathbf{y}}(x_{i+1})| + \cdots + |f_{\mathbf{y}}(x_{j-1}) - f_{\mathbf{y}}(x')|$$
$$\leq 2Q|x - x_i|^{\tau} + 2Q|x_i - x_{i+1}|^{\tau} + \cdots + 2Q|x_{j-1} - x'|^{\tau}$$
$$\leq 2nQ|x - x'|^{\tau}$$
$$\leq K|x - x'|^{\tau},$$
where $K = 2nQ$. □

Theorem 5 shows that if there are constants $A$ and $\rho$ such that $|y_i| \leq A$, $L_{q_i} \leq \rho$, then the fractal interpolation function $f_{\mathbf{y}}$ is a Hölder continuous function of which Hölder coefficient $K$ and Hölder exponent $\tau$ are independent of $y_i$ and $q_i$.

Now, we discuss with Hölder continuity of the fuzzy-valued fractal interpolation function $f \in C^*(I,\ \mathbf{R_F})$.

**Theorem 6.** *The fuzzy-valued fractal interpolation function $f$ is Hölder continuous, that is, there exist a constant $H_f$ and a constant $0 < \tau \leq 1$ such that*

$$\forall x, x' \in I, \ d_\infty(f(x), f(x')) \leq H_f \, |x-x'|^\tau.$$

**Proof**. From Theorem 4, for any $0 \leq \lambda \leq 1$, $(f(\cdot))^-(\lambda)$ and $(f(\cdot))^+(\lambda)$ are the fractal interpolation functions $f_{\mathbf{u}^-(\lambda)}$ and $f_{\mathbf{u}^+(\lambda)}$ interpolating the data sets $P_\lambda^-$ and $P_\lambda^+$, respectively, where

$$\mathbf{u}^-(\lambda) = (u_0^-(\lambda), u_1^-(\lambda), \cdots, u_n^-(\lambda)) \text{ and } \mathbf{u}^+(\lambda) = (u_0^+(\lambda), u_1^+(\lambda), \cdots, u_n^+(\lambda)).$$

Moreover,

$$\min_{i=0,1,\cdots,n}\{u_i^-(0)\} \leq u_i^-(\lambda) \leq \max_{i=0,1,\cdots,n}\{u_i^-(0)\}, \quad \min_{i=0,1,\cdots,n}\{u_i^+(0)\} \leq u_i^+(\lambda) \leq \max_{i=0,1,\cdots,n}\{u_i^+(0)\}$$

and since $q_i$ are Lipschitz functions satisfying $L_{q_i} \leq \rho$, $(q_i(\cdot))^\pm(\lambda)$ are also Lipschitz functions satisfying

$$L_{(q_i(\cdot))^\pm(\lambda)} \leq \rho.$$

Hence, by Theorem 5, there exist constants $K^-$, $K^+$ and $0 < \tau \leq 1$ independent of $\lambda (\in [0,1])$ and $\mathbf{u}^-(\lambda), \mathbf{u}^+(\lambda)$ such that for any $x, x' \in I$, we have

$$|f_{\mathbf{u}^-(\lambda)}(x) - f_{\mathbf{u}^-(\lambda)}(x')| \leq K^- |x-x'|^\tau,$$

$$|f_{\mathbf{u}^+(\lambda)}(x) - f_{\mathbf{u}^+(\lambda)}(x')| \leq K^+ |x-x'|^\tau.$$

Then we have

$$\begin{aligned}
d_\infty(f(x), f(x')) &= \sup_{\lambda \in [0,1]} \max\{|(f(x))^-(\lambda) - (f(x'))^-(\lambda)|, |(f(x))^+(\lambda) - (f(x'))^+(\lambda)|\} \\
&= \sup_{\lambda \in [0,1]} \max\{|f_{\mathbf{u}^-(\lambda)}(x) - f_{\mathbf{u}^-(\lambda)}(x')|, |f_{\mathbf{u}^+(\lambda)}(x) - f_{\mathbf{u}^+(\lambda)}(x')|\} \\
&\leq \max\{K^- |x-x'|^\tau, K^+ |x-x'|^\tau\} \\
&\leq \max\{K^-, K^+\} |x-x'|^\tau \\
&= H_f |x-x'|^\tau,
\end{aligned}$$

where $H_f = \max\{K^-, K^+\}$. □

## 6. Conclusion

There are many objects with irregularity in nature. Since the data obtained from them can be perturbed by some factors, they have not only the irregularity but also the uncertainty. To better model them, firstly, we constructed a fuzzy-valued fractal interpolation function that interpolates a data set of fuzzy numbers using an iterated function system. Next, we proved that it is Hölder continuous. In the future, on the basis of the results, we are going to study the stability and derivatives and integrals of the fuzzy-valued fractal interpolation functions.

### References

1. M. F. Barnsley Fractal functions and interpolation, *Constr.Approx.2* (1986) 303-329
2. Z. G. Feng, H. P. Xie, On stability of fractal interpolation, *Fractals* 6(3) (1998) 269-273.
3. J. Z. Xiao, X. H. Zhu and P. P. Jin, Iterated function systems and attractors in the KM fuzzy metric spaces, *Fuzzy Sets and Systems* 267 (2015) 100–116.
4. C. H. Yun, Hidden variable recurrent fractal interpolation function with four function contractivity factors, *Fractals* 27(7) (2019) 1950113 1-13.

5. C. H. Yun, M. G. Ri, Box-counting dimension and analytic properties of hidden variable fractal interpolation functions with function contractivity factors, *Chaos Solitons Fractals* 134 (2020) 109700.
6. L. A. Zadeh, Fuzzy sets, *Inform. Control* 8 (1965) 338–353.
7. R. Lowen, A fuzzy lagrange interpolation theorem, *Fuzzy sets and systems* 34 (1990) 33-38.
8. O. Kaleva, Interpolation of fuzzy data, *Fuzzy Sets and Systems* 61 (1994) 63-70.
9. H. S. Goghary, S. Abbasbandy, Interpolation of fuzzy data by Hermite polynomial, *International Journal of Computer Mathematics* 82(12) (2005) 1541–1545.
10. S. Abbasbandy, R. Ezzati and H. Behforooz, Interpolation of fuzzy data by using fuzzy splines, *International Journal of Uncertainty, Fuzziness and Knowledge-Based Systems* 16(1) (2008) 107–115.
11. J. Andres, M. Rypka, Fuzzy fractals and hyperfractals, *Fuzzy Sets and Systems* http://dx.doi.org/10.1016/j.fss.2016.01.008 (2016) 1-17.
12. C. A. Cabrelli, B. Forte, U. M. Molter and E. R. Vrscay Iterated fuzzy set systems: a new approach to the inverse problem for fractals and other sets, *J. Math. Anal. Appl.* 171 (1992) 79–100.
13. D. Easwaramoorthy, R. Uthayakumar, Analysis on fractals in fuzzy metric spaces, *Fractals* 19 (3) (2011) 379–386.
14. K. K. Majumdar, Fuzzy fractals and fuzzy turbulence, *IEEE Trans. Syst. Man Cybern., Part B, Cybern.* 34 (1) (2004) 746–752.
15. R. Uthayakumar, D. Easwaramoorthy, Hutchinson–Barnsley operator in fuzzy metric spaces, *Int. J. Eng. Nat. Sci.* 5(4) (2011) 203–207.
16. X. Xiao, Z. Li and S. Yan, Fitting of fuzzy fractal interpolation for uncertain data, *CCIS* 86 (2011) 78-84.
17. X. Xiao, Z. Li and W. Gong, Fuzzy fractal interpolation surface and its applications, *Advanced Materials Research* 542-543 (2012) 1141-1144.
18. J. X. Fang, Q. Y. Xue, Some properties of the space of fuzzy-valued continuous functions on a compact set, *Fuzzy Sets and Systems* 160 (2009) 1620–1631.